\documentclass[11pt,a4paper]{article}

\usepackage{amssymb,amsmath,amsthm,indentfirst}

\newtheorem{thm}{Theorem}[section]
\newtheorem{prop}[thm]{Proposition}
\newtheorem{lem}[thm]{Lemma}

\theoremstyle{definition}
\newtheorem{defi}[thm]{Definition}
\newtheorem{rem}[thm]{Remark}

\newcommand{\f}[1]{\tilde{f}_{#1}}
\newcommand{\e}[1]{\tilde{e}_{#1}}

\DeclareMathOperator{\Hom}{Hom}
\DeclareMathOperator{\Ker}{Ker}
\DeclareMathOperator{\Ima}{Im}
\DeclareMathOperator{\wt}{wt}
\DeclareMathOperator{\jw}{Jw}

\renewcommand\theenumi{\roman{enumi}}
\renewcommand\labelenumi{{\rm(\theenumi)}}

\title{\textbf{A generalization of adjoint crystals\\ for the quantized affine algebras\\ of type $A \sb{n} \sp{(1)}$, $C \sb{n} \sp{(1)}$ and $D \sb{n+1} \sp{(2)}$}}

\author{Ryosuke Kodera}

\date{}

\makeatletter
\let\@old@@maketitle=\@maketitle
\def\@maketitle{%
\footnotetext{%
\hspace*{-1em}\hspace*{-\footnotesep}%
Graduate School of Mathematical Sciences, The University of Tokyo, 3-8-1 Komaba, Meguro, Tokyo 153-8914, Japan.\\
E-mail address: kryosuke@ms.u-tokyo.ac.jp\\
2000 \textit{Mathematics Subject Classification.} 17B37.}
\@old@@maketitle
}
\makeatother

\begin{document}
\maketitle

\begin{abstract}
We propose to generalize Benkart-Frenkel-Kang-Lee's adjoint crystals and describe their crystal structure for type $A \sb{n} \sp{(1)}$, $C \sb{n} \sp{(1)}$ and $D \sb{n+1} \sp{(2)}$. 
\end{abstract}

\section{Introduction}
Let $U_q'(\mathfrak{g})$ be a quantized affine algebra without the degree operator.
Let $\mathfrak{g}_0$ be the underlying simple Lie algebra of finite type contained in the affine Lie algebra $\mathfrak{g}$ and $U_q(\mathfrak{g}_0)$ the corresponding subalgebra of $U_q'(\mathfrak{g})$.
Let $i$ be a node of the Dynkin diagram of $\mathfrak{g}_0$ and $l$ a nonnegative integer.
We denote by $W \sp{i,l}$ the Kirillov-Reshetikhin module for $i$ and $l$, which is a finite-dimensional irreducible module over $U_q'(\mathfrak{g})$.
They are important objects from the viewpoint of the crystal base theory as it is conjectured that all Kirillov-Reshetikhin modules have crystal bases.

We assume that $\mathfrak{g}$ is not of type $A \sb{2n} \sp{(2)}$ for simplicity.
We denote by $V(\lambda)$ the finite-dimensional irreducible module with highest weight $\lambda$ over $U_q(\mathfrak{g}_0)$ and $B(\lambda)$ its crystal base.
In \cite{bfkl}, they study a $U_q'(\mathfrak{g})$-module $V$ with the following properties.
\begin{itemize}
\item $V$ decomposes as $V(\theta) \oplus V(0)$ as a $U_q(\mathfrak{g}_0)$-module. 
\item $V$ has a crystal base.
\end{itemize}
Here $\theta$ is given by
\[\theta=
\begin{cases}
\varpi_1 + \varpi_n         & \text{for type } A \sb{n} \sp{(1)},\\
2\varpi_1                   & \text{for type } C \sb{n} \sp{(1)},\\
\varpi_{i_0}                & \text{for other types,}
\end{cases}\]     
where $\varpi_i$ denotes the $i$-th fundamental weight of $\mathfrak{g}_0$ and $i_0$ denotes the node of the Dynkin diagram of $\mathfrak{g}$ connected to the special node $0$.
If $\mathfrak{g}$ is an untwisted affine algebra, then $\theta$ is the highest root of $\mathfrak{g}_0$.
Otherwise, namely if $\mathfrak{g}$ is a twisted affine algebra, then $\theta$ is the highest short root of $\mathfrak{g}_0$.
The structure of the crystal base of $V$ is determined by defining the action of the affine Kashiwara operator $\f{0}$ on $B(\theta) \oplus B(0)$ correctly for all quantized affine algebras in a uniform manner.
We call them adjoint crystals here.

As a natural generalization of \cite{bfkl}, we consider a family of $U_q'(\mathfrak{g})$-modules $\{ V_l \} \sb{l \in \mathbb{Z}_{\geq 0}}$ which has the following properties.
\begin{itemize}
\item $V_l$ decomposes as $\bigoplus \sb{k=0} \sp{l} V(k\theta)$ as a $U_q(\mathfrak{g}_0)$-module. 
\item $V_l$ has a crystal base.
\end{itemize}
In fact, $V \sb{l}$ is given as follows:
\[V_l=
\begin{cases}
W \sp{1,l} \otimes W \sp{n,l} & \text{for type } A \sb{n} \sp{(1)},\\
W \sp{1,2l}                   & \text{for type } C \sb{n} \sp{(1)},\\
W \sp{i\sb{0},l}              & \text{for other types.}
\end{cases}\]     
The decomposition of $V_l$ as a $U_q(\mathfrak{g}_0)$-module follows from results of \cite{chari} for untwisted cases and \cite{her} for twisted cases.
Moreover we can easily see that $V_l$ satisfies the sufficient condition for the existence of a crystal base stated in \cite[Proposition~3.4.5]{kmn}.
We denote by $B_l$ the crystal base of $V_l$.
Remark that the crystal $B_l$ has been studied for some cases. 
For example, the case of type $C \sb{n} \sp{(1)}$ appears in \cite{kkm}, $D \sb{n} \sp{(1)}$ in \cite{schstern}, $D \sb{n+1} \sp{(2)}$ in \cite{kmn}, $G \sb{2} \sp{(1)}$ in \cite{yamane} and $D \sb{4} \sp{(3)}$ in \cite{yamada}.
Inspired by \cite{bfkl}, we expect that (i) the crystal graph of $B_{l-1}$ is a full subgraph of that of $B_l$ and (ii) there exists a simple rule for extending the action of $\f{0}$ on $B_{l-1}$ to that on $B_l$.
We consider in the present paper the case of $A \sb{n} \sp{(1)}$, $C \sb{n} \sp{(1)}$ and $D \sb{n+1} \sp{(2)}$ to show that the above expectations (i) and (ii) are true in these cases.  

Let us explain results of this paper for type~$A \sb{n} \sp{(1)}$.
Similar assertions also hold for type~$C \sb{n} \sp{(1)}$ and $D \sb{n+1} \sp{(2)}$.
By the decomposition of $B_{l}$ as a $U_q(\mathfrak{g}_0)$-crystal, we regard $B(k\theta)$ as a subset of $B_l$.
We define $n+1$ maps $\Theta_1, \dots, \Theta_{n+1}$ from $B_{l-1}$ to $B_l$.
Our results are summarized as follows.

\begin{thm}
\begin{enumerate}
\item The crystal graph of $B_{l-1}$ is regarded as a full subgraph of that of $B_l$ via $\Theta_1$.

\item For $j=2, \dots, n+1$, the map $\Theta_j$ and the Kashiwara operator $\f{0}$ commute with each other.

\item We have 
\[\bigcup_{j=2}^{n+1} \Ima\Theta_j = \bigsqcup_{k=0}^l \, \{b \in B(k\theta) \subset B_l \mid \wt b \in \wt B((k-1)\theta)\}\]
and the space $B_l \setminus \bigcup_{j=2}^{n+1} \Ima\Theta_j$ is weight multiplicity free.

\item Let $b \in B_l \setminus \bigcup_{j=2}^{n+1} \Ima\Theta_j$.
Then the element $\f{0}b$ is uniquely determined by its weight.
\end{enumerate}
\end{thm}
 
\subsection*{Acknowledgments}
The author would like to thank his advisor Yoshihisa Saito for his guidance and valuable comments.
He also would like to thank Masato Okado for answering the author's questions about Kirillov-Reshetikhin crystals, and Noriyuki Abe for providing a computer program for drawing crystal graphs.  

\section{Preliminaries}\label{sec:pre}

\subsection{Quantized universal enveloping algebras and crystal bases}
We shall review on quantized universal enveloping algebras and the crystal base theory based on \cite{kasc}.
We also refer \cite{hongkang}.

Suppose that the following data are given:
\begin{quote}
$P$: a free $\mathbb{Z}$-module,

$I$: an index set,

$\Pi = \{\alpha_i\mid i\in I\}\subset P$,

$\Pi^\vee = \{h_i\mid i\in I\}\subset P^\vee=\Hom_{\mathbb{Z}}(P,\mathbb{Z})$,

$(\,\cdotp,\cdotp)$: a $\mathbb{Q}$-valued symmetric bilinear form on $P$.
\end{quote}
These data are supposed to satisfy the following conditions:
\begin{quote}
$(\alpha_i,\alpha_i)>0$\quad for any $i\in I$,

$(\alpha_i,\alpha_j)\leq 0$\quad for any $i,j\in I$ with $i\neq j$,

$\langle h_i,\lambda\rangle=\dfrac{2(\alpha_i,\lambda)}{(\alpha_i,\alpha_i)}$\quad for any $i\in I$ and $\lambda\in P$.
\end{quote}
We do not assume $\Pi$ and $\Pi \sp{\vee}$ to be linearly independent sets in general.
Note that $(\langle h_i,\alpha_j\rangle)_{i,j\in I}$ is a symmetrizable generalized Cartan matrix.
Let $\mathfrak{g}$ be the associated Kac-Moody Lie algebra.

Let $\gamma$ be the minimal positive integer such that $(\alpha_i,\alpha_i)/2\in \gamma^{-1}\mathbb{Z}$ for any $i\in I$.
Let $q$ be an indeterminate over $\mathbb{Q}$ and put $q_s=q^{1/\gamma}$.
We use the notation:
\[q_i = q^{\frac{(\alpha_i, \alpha_i)}{2}}, \, [k]_i = \dfrac{q_i^k-q_i^{-k}}{q_i-q_i^{-1}}, \, [k]_i \,! = \prod_{r=1}^k \, [r]_i.\]

\begin{defi}
The quantized universal enveloping algebra $U_q(\mathfrak{g})$ associated with $(P,P^\vee,\Pi,\Pi^\vee)$ is the unital associative algebra over the rational fraction field $\mathbb{Q}(q_s)$ generated by $e_i$ and $f_i$ for $i\in I$ and $q^h$ for $h\in \gamma^{-1}P^{\vee}$ with the following defining relations:
\begin{enumerate}
\item $q^0=1,\ q^h q^{h'}=q^{h+h'}$\quad for $h,h'\in \gamma^{-1}P^{\vee}$,

\item $q^h e_i q^{-h}=q^{\langle h,\alpha_i\rangle} e_i$\quad for $i\in I$ and $h\in \gamma^{-1}P^{\vee}$,

\item $q^h f_i q^{-h}=q^{-\langle h,\alpha_i\rangle} f_i$\quad for $i\in I$ and $h\in \gamma^{-1}P^{\vee}$,

\item $[e_i,f_j]=\delta_{ij}\dfrac{t_i -{t_i}^{-1}}{q_i-{q_i}^{-1}}$\quad for $i,j\in I$,

\item $\displaystyle\sum_{k=0}^{1-\langle h_i,\alpha_j\rangle} (-1)^k {e_i}^{(k)} e_j {e_i}^{(1-\langle h_i,\alpha_j\rangle-k)}=0$\quad for $i,j\in I$ with $i\neq j$,

\item $\displaystyle\sum_{k=0}^{1-\langle h_i,\alpha_j\rangle} (-1)^k {f_i}^{(k)} f_j {f_i}^{(1-\langle h_i,\alpha_j\rangle-k)}=0$\quad for $i,j\in I$ with $i\neq j$.
\end{enumerate}
Here $t_i = q^{\frac{(\alpha_i, \alpha_i)}{2}h_i}$, $e_i^{(k)} = \dfrac{e_i^k}{[k]_i\,!}$ and $f_i^{(k)} = \dfrac{f_i^k}{[k]_i\,!}$.
\end{defi}

We call $P$, $P^\vee$, $\alpha_i$ and $h_i$ the weight lattice, the coweight lattice, a simple root and a simple coroot, respectively.
We denote by $P^+$ the set of dominant weights, that is, $P^+ = \{ \lambda \in P \mid \langle h_i, \lambda \rangle \geq 0 \text{ for any } i \in I\}$.
We define the root lattice $Q$ and its positive part $Q^+$ by $Q=\sum_{i\in I}\mathbb{Z}  \alpha_i$ and $Q^+=\sum_{i\in I}\mathbb{Z}_{\geq 0} \alpha_i$.
For $\lambda,\mu \in P$, we say $\lambda \geq \mu$ if $\lambda - \mu\in Q^+$.  

Let $M$ be an integrable $U_q(\mathfrak{g})$-module and $M=\bigoplus_{\lambda\in P}M_\lambda$ its weight space decomposition.
According to the representation theory of $U_q(\mathfrak{sl}_2)$, $M_\lambda$ decomposes as
\[M_\lambda=\bigoplus_{n\geq 0}f_i^{(n)}(\Ker e_i \cap M_{\lambda+n\alpha_i})\]
for each $i$.
We define endomorphisms $\e{i}$ and $\f{i}$ of $M$ by
\[ \e{i}(f_i^{(n)}u)=f_i^{(n-1)}u \]
and
\[ \f{i}(f_i^{(n)}u)=f_i^{(n+1)}u \]
for $u \in \Ker e_i \cap M_{\lambda+n\alpha_i}$.
They are called Kashiwara operators or modified root operators.

Let $A$ be the subring of $\mathbb{Q}(q_s)$ that consists of rational fractions without a pole at $q_s=0$.   

\begin{defi}
A pair $(L,B)$ is called a crystal base of an integrable $U_q(\mathfrak{g})$-module $M$ if it satisfies the following conditions: 
\begin{enumerate}
\item $L$ is a free $A$-submodule of $M$ such that $L\otimes_A \mathbb{Q}(q_s)\simeq M$,

\item $B$ is a $\mathbb{Q}$-basis of $L/q_s L$,

\item $\e{i}L \subset L$ and $\f{i}L \subset L$ for any $i\in I$, hence $\e{i}$ and $\f{i}$ act on $L/q_s L$, 

\item $\e{i}B \subset B\sqcup \{0\}$ and $\f{i}B \subset B \sqcup \{0\}$ for any $i\in I$, 

\item $L=\bigoplus_{\lambda} L_\lambda$ where $L_\lambda=L\cap M_\lambda$,

\item $B=\bigsqcup_{\lambda} B_\lambda$ where $B_\lambda=B\cap (L_\lambda/q_s L_\lambda)$,

\item for $b,b'\in B$, $b'=\f{i}b$ if and only if $b=\e{i}b'$.
\end{enumerate}
\end{defi}

We often regard $B$ as a crystal base rather than $(L,B)$.
For a crystal base $B$, we define a colored oriented graph called the crystal graph as follows.
The vertices of the graph are elements of $B$.
For $b,b'\in B$, draw an arrow labeled by $i$ from $b$ to $b'$ if $b'=\f{i}b$.

Let $\wt\colon B\to P$ be the map such that $\wt b =\lambda$ for $b\in B_\lambda$.
For $b\in B$ and $i\in I$, we set 
\[\varepsilon_i(b)=\max\{n\mid \e{i}^n b \neq 0\}\quad \text{and}\quad \varphi_i(b)=\max\{n\mid \f{i}^n b \neq 0\}.\]

Let $M_1$ and $M_2$ be integrable $U_q(\mathfrak{g})$-modules with crystal bases $(L_1,B_1)$ and $(L_2,B_2)$.
We denote by $B_1 \oplus B_2$ a direct sum $B_1 \sqcup B_2$.
Then $(L_1 \oplus_A L_2, B_1 \oplus B_2)$ gives a crystal base of $M_1 \oplus M_2$.
We denote by $B_1\otimes B_2$ a direct product $B_1\times B_2$, which is a $\mathbb{Q}$-basis of $L_1/q_s L_1 \otimes_{\mathbb{Q}} L_2/q_s L_2 \simeq (L_1 \otimes_A L_2)/q_s (L_1 \otimes_A L_2)$.
Then $(L_1\otimes_A L_2,B_1\otimes B_2)$ is a crystal base of $M_1\otimes M_2$.
Kashiwara operators act by
\[ \e{i}(b_1 \otimes b_2)=
\begin{cases}
(\e{i}b_1) \otimes b_2 & \text{if} \ \varphi_i(b_1) \geq \varepsilon_i(b_2),\\
b_1 \otimes (\e{i}b_2) & \text{if} \ \varphi_i(b_1) < \varepsilon_i(b_2),
\end{cases}\]
\[ \f{i}(b_1 \otimes b_2)=
\begin{cases}
(\f{i}b_1) \otimes b_2 & \text{if} \ \varphi_i(b_1) > \varepsilon_i(b_2),\\
b_1 \otimes (\f{i}b_2) & \text{if} \ \varphi_i(b_1) \leq \varepsilon_i(b_2).
\end{cases}\]

The notion of crystals is a combinatorial generalization of crystal bases.
We shall not review the theory of abstract crystals here since all crystals appearing in this article come from crystal bases.
It can be found in \cite[Section~4.5]{hongkang}.

A bijection between two crystal bases is called an isomorphism of crystals if it commutes with Kashiwara operators and preserves weights.

\subsection{Quantized affine algebras}

Now we shall turn to affine situations.
See \cite[Chapter~6-8]{kac} for a reference on affine Lie algebras.
Let $\mathfrak{g}$ be an affine Lie algebra with a Cartan matrix indexed by $I$.
We assume that $\mathfrak{g}$ is not of type $A_{2n}^{(2)}$ for simplicity.
Let $\delta$ be the generator of imaginary roots, $c$ the canonical central element and $d$ the degree operator. 
Choose $0 \in I$ such that $\delta - \alpha_0 \in \sum \sb{i \in I, i \neq 0} \mathbb{Z} \alpha_i$ and set $I_0 = I \setminus \{0\}$.
Let $\mathfrak{g}_0$ be the underlying simple Lie algebra of finite type.
Let $\Lambda_i$ be the fundamental weight for $i \in I$ and define the weight lattice $P$ by 
\[P=\bigoplus_{i\in I} \mathbb{Z}\Lambda_i \oplus \mathbb{Z} \delta.\]
We normalize the symmetric invariant bilinear form on $P$ so that $(\delta,\lambda)=\langle c,\lambda \rangle$ for any $\lambda \in P$.
Set $P_{\text{cl}}=P / \mathbb{Z}\delta$ and denote by the same letter the image of $\Lambda_i$ for each $i$. 
Then we have
\[P\sb{\text{cl}} = \bigoplus_{i\in I} \mathbb{Z}\Lambda_i.\]
We define the level zero fundamental weight $\varpi_i$ for $i\in I_0$, which is an element of $P_{\text{cl}}$, by $\varpi_i = \Lambda_i - \langle c, \Lambda_i \rangle \Lambda_0$.
We set $P_0= P_{\text{cl}} / \mathbb{Z}\Lambda_0$ and denote by the same letter the image of $\varpi_i$ for each $i$.
Then we have
\[P_0 = \bigoplus_{i\in I_0} \mathbb{Z}\varpi_i.\]
We identify $P_0$ with the weight lattice of $\mathfrak{g}_0$ and $\varpi_i$ its $i$-th fundamental weight.
Then quantized algebras $U_q(\mathfrak{g})$, $U_q'(\mathfrak{g})$ and $U_q(\mathfrak{g}_0)$ are defined associated with weight lattices $P$, $P_{\text{cl}}$ and $P_0$.  
Note that simple roots in $P_{\text{cl}}$ are linearly dependent while they are linearly independent in $P$ and $P_0$.
However we use the same letters for simple roots in $P \sb{\text{cl}}$ and $P_0$.
For the coweight lattices, we have
\[P_0^\vee = \bigoplus \sb{i\in I_0} \mathbb{Z} h_i
\, \subset \, P_{\text{cl}}^\vee = \bigoplus \sb{i\in I} \mathbb{Z} h_i
\, \subset \, P^\vee = \bigoplus \sb{i\in I} \mathbb{Z} h_i \oplus \mathbb{Z} d. \]
Hence $U_q'(\mathfrak{g})$ is the subalgebra of $U_q(\mathfrak{g})$ generated by $e_i$ and $f_i$ for $i\in I$ and $q^h$ for $h\in \gamma^{-1}P_{\text{cl}}^\vee$ and $U_q(\mathfrak{g}_0)$ is the subalgebra of $U_q'(\mathfrak{g})$ generated by $e_i$ and $f_i$ for $i\in I_0$ and $q^h$ for $h\in \gamma^{-1} P_0^\vee$.
When we need to clarify in which weight lattice we work, write a $U_q(\mathfrak{g}_0)$-crystal, a $U_q(\mathfrak{g}_0)$-weight, etc.  

Since the definition of Kirillov-Reshetikhin modules is not used in this paper, we omit it and refer \cite[Section~3]{okadosch}.
They are finite-dimensional irreducible $U_q'(\mathfrak{g})$-modules.
We denote by $W \sp{i,l}$ the Kirillov-Reshetikhin module for $i \in I_0$ and a nonnegative integer $l$.
For the quantized affine algebras of nonexceptional types, it is proved that any Kirillov-Reshetikhin module has a crystal base in \cite{kmn} and \cite{okadosch}.
They are called Kirillov-Reshetikhin crystals and denoted by $B^{i,l}$.

\subsection{Some lemmas on weights}\label{subsec:lemma}

Assume that $\mathfrak{g}_0$ is either of type $A_n$, $C_n$ or $B_n$ in this subsection.
For $\lambda \in P_0^+$, we denote by $B(\lambda)$ the crystal base of the finite-dimensional irreducible $U_q(\mathfrak{g}_0)$-module with highest weight $\lambda$.
If $\mathfrak{g}_0$ is of type $A_n$ or $C_n$, we denote by $\theta$ the highest root of $\mathfrak{g}_0$.
If $\mathfrak{g}_0$ is of type $B_n$, denote by $\theta$ the highest short root of $\mathfrak{g}_0$.
We shall prove some lemmas on weights of $B(k\theta)$.

It is well known that
\[\wt B(\lambda) = W \cdot \{\nu \in P_0^+ \mid \nu \leq \lambda\}\]
for any $\lambda \in P_0^+$.
Here $W$ is the Weyl group of $\mathfrak{g}_0$.
(See e.g.\ \cite[21.3~Proposition]{hum}.)
Hence the following lemma is immediate.
\begin{lem}\label{lem:inclusion}
Let $k$ and $k'$ be nonnegative integers with $k' \leq k$.
Then $\wt B(k\theta)$ contains $\wt B(k'\theta)$. 
\end{lem}

We use the following standard notation for the finite root systems of type $A_n$, $C_n$ and $B_n$.
In the case of type $A_n$, simple roots and fundamental weights are defined by
\begin{align*}
\alpha_i &= \epsilon_i - \epsilon_{i+1} \ \text{ for } i \in I_0,\\
\varpi_i &= \sum_{j=1}^{i} \epsilon_j \ \text{ for } i \in I_0.
\end{align*}
Here we define the element $\epsilon_j \in P_0$ for $j = 1, \dots, n+1$ by
\[\langle h_i, \epsilon_j \rangle =
\begin{cases}
1  & \text{if } j = i,\\
-1 & \text{if } j = i+1,\\
0  & \text{otherwise}.
\end{cases}
\]
Note that $\sum_{j=1}^{n+1} \epsilon_j = 0$.
The highest root $\theta$ is given by
\[\theta = \alpha_1 + \cdots + \alpha_n = \epsilon_1 - \epsilon_{n+1} = \varpi_1 + \varpi_n.\]
In the case of type $C_n$,
\begin{align*}
\alpha_i &= \epsilon_i - \epsilon_{i+1} \ \text{ for } i = 1, \dots, n-1,\\
\alpha_n &= 2\epsilon_n,\\
\varpi_i &= \sum_{j=1}^{i} \epsilon_j \ \text{ for } i \in I_0,
\end{align*}
where $\epsilon_j \in P_0$ for $j = 1, \dots, n$ is defined by
\[\langle h_i, \epsilon_j \rangle =
\begin{cases}
1  & \text{if } j = i,\\
-1 & \text{if } j = i+1,\\
0  & \text{otherwise}.
\end{cases}
\]
The highest root $\theta$ is given by
\[\theta = 2\alpha_1 + \cdots + 2\alpha_{n-1} + \alpha_n = 2\epsilon_1 = 2\varpi_1.\]
In the case of type $B_n$,
\begin{align*}
\alpha_i &= \epsilon_i - \epsilon_{i+1} \ \text{ for } i = 1, \dots, n-1,\\
\alpha_n &= \epsilon_n,\\
\varpi_i &= \sum_{j=1}^{i} \epsilon_j \ \text{ for } i = 1, \dots, n-1,\\
\varpi_n &= \dfrac{1}{2}\sum_{j=1}^n\epsilon_j,
\end{align*}
where $\epsilon_j \in P_0$ for $j = 1, \dots, n$ is defined by
\[\langle h_i, \epsilon_j \rangle =
\begin{cases}
1  & \text{if } j = i,\\
-1 & \text{if } j = i+1,\\
0  & \text{otherwise}
\end{cases}
\]
for $i=1, \dots, n-1$ and 
\[\langle h_n, \epsilon_j \rangle =
\begin{cases}
2  & \text{if } j = n,\\
0  & \text{otherwise}.
\end{cases}
\]
The highest short root $\theta$ is given by
\[\theta = \alpha_1 + \cdots + \alpha_n = \epsilon_1 = \varpi_1.\]

Let $\mu \in P_0$.
We denote by $m_j(\mu)$ the coefficient of $\epsilon_j$ in $\mu$.
For type $A_n$, we normalize them so that $\sum_{j=1}^{n+1} m_j(\mu) = 0$.
For type $A_n$, set 
\[J(\mu) = \{j \mid m_j(\mu) > 0 \}.\]
For type $C_n$ and $B_n$, set
\[|\mu| = \sum_{j=1}^{n} |m_j(\mu)|.\]

\begin{lem}\label{lem:sum}
Let $k$ be a nonnegative integer.
\begin{enumerate}
\item Let $\mathfrak{g}_0$ be of type $A_n$.
For $\mu \in P_0$, $\mu \in \wt B(k\theta)$ if and only if $\sum_{j \in J(\mu)} m_j(\mu) \leq k$.
In particular, $\mu \in \wt B(k\theta) \setminus \wt B((k-1)\theta)$ if and only if $\sum_{j \in J(\mu)}  m_j(\mu) = k$.

\item Let $\mathfrak{g}_0$ be of type $C_n$.
For $\mu \in P_0$, $\mu \in \wt B(k\theta)$ if and only if $2k - |\mu| \in 2\mathbb{Z}_{\geq 0}$.
In particular, $\mu \in \wt B(k\theta) \setminus \wt B((k-1)\theta)$ if and only if $|\mu| = 2k$.

\item Let $\mathfrak{g}_0$ be of type $B_n$.
For $\mu \in P_0$, $\mu \in \wt B(k\theta)$ if and only if $|\mu| \leq k$.
In particular, $\mu \in \wt B(k\theta) \setminus \wt B((k-1)\theta)$ if and only if $|\mu| = k$.
\end{enumerate}
\end{lem}

\begin{proof}
We prove (i).
For $\mu \in P_0$, the coefficient of $\alpha_i$ in $\mu$ is given by $\sum_{j=1}^{i}m_j (\mu)$.
Therefore $\mu \in \wt B(k\theta)$ if and only if there exists an element $\tau$ of the symmetric group of degree $n+1$ satisfying
\[m_{\tau(1)}(\mu) \geq m_{\tau(2)}(\mu) \geq \dots \geq m_{\tau(n+1)}(\mu)\]
and
\[ k \geq \sum_{j=1}^{i} m_{\tau(j)}(\mu)\]
for $i = 1, \dots, n$.
Since
\[\sum_{j\in J(\mu)} m_j(\mu) = \max \{\sum_{j=1}^i m_{\tau(j)}(\mu) \mid 1 \leq i \leq n\},\]
$\mu \in \wt B(k\theta)$ if and only if
\[k \geq \sum_{j\in J(\mu)} m_j(\mu).\]

We prove (ii).
The coefficient of $\alpha_i$ in $\mu$ is given by $\sum_{j=1}^{i} m_j(\mu)$ for $i \neq n$ and $(\sum_{j=1}^{n} m_j(\mu))/2$ for $i=n$.
Assume that $\mu$ is dominant, that is, 
\[m_1(\mu) \geq m_2(\mu) \geq \dots \geq m_n(\mu) \geq 0.\] 
Then $\mu \in \wt B(k\theta)$ if and only if
\[2k \geq \sum_{j=1}^{i} m_j(\mu)\]
for $i = 1, \dots, n-1$ and  
\[k - \dfrac{1}{2}\sum_{j=1}^{n} m_j(\mu) \in \mathbb{Z}_{\geq 0}.\]
This condition is equivalent to 
\[2k - \sum_{j=1}^{n} m_j(\mu) \in 2\mathbb{Z}_{\geq 0}.\]
For general $\mu$, $\mu$ is $W$-conjugate to a dominant weight with the above condition if and only if $\mu$ satisfies
\[2k - \sum_{j=1}^{n} |m_j(\mu)| \in 2\mathbb{Z}_{\geq 0}.\]   

We prove (iii).
The coefficient of $\alpha_i$ in $\mu$ is given by $\sum_{j=1}^{i} m_j(\mu)$.
Assume that $\mu$ is dominant, that is, 
\[m_1(\mu) \geq m_2(\mu) \geq \dots \geq m_n(\mu) \geq 0.\] 
Then $\mu \in \wt B(k\theta)$ if and only if
\[k \geq \sum_{j=1}^i m_j(\mu)\]
for $i = 1, \dots, n$.
This condition is equivalent to 
\[k \geq \sum_{j=1}^n m_j(\mu).\]
For general $\mu$, $\mu$ is $W$-conjugate to a dominant weight with the above condition if and only if $\mu$ satisfies
\[k \geq \sum_{j=1}^{n} |m_j(\mu)|.\]   
\end{proof}

Let $\mu \in \wt B(k\theta) \setminus \wt B((k-1)\theta)$. 
Now we investigate the behavior of $\mu + \theta$.
Let $\mathfrak{g}_0$ be of type $A_n$.
There are the following four cases:
\begin{enumerate}
\renewcommand\theenumi{\alph{enumi}}
\renewcommand\labelenumi{\rm{(\theenumi)}}
\item $m_1(\mu) \geq 0$ and $m_{n+1}(\mu) \leq 0$,
\item $m_1(\mu) \geq 0$ and $m_{n+1}(\mu) > 0$,
\item $m_1(\mu) < 0$ and $m_{n+1}(\mu) \leq 0$,
\item $m_1(\mu) < 0$ and $m_{n+1}(\mu) > 0$.
\end{enumerate}

\begin{lem}\label{lem:wt}
Let $\mathfrak{g}_0$ be of type $A_n$ and assume $\mu \in \wt B(k\theta) \setminus \wt B((k-1)\theta)$.
Then we have
\begin{enumerate}
\item $\mu + \theta \in \wt B((k+1)\theta) \setminus \wt B(k\theta)$ if and only if $\mu$ satisfies {\rm(a)}.

\item $\mu + \theta \in \wt B(k\theta) \setminus \wt B((k-1)\theta)$ if and only if $\mu$ satisfies {\rm(b)} or {\rm(c)}.

\item $\mu + \theta \in \wt B((k-1)\theta) \setminus \wt B((k-2)\theta)$ if and only if $\mu$ satisfies {\rm(d)}. 
\end{enumerate}					
\end{lem}

\begin{proof}
Since $\theta =\epsilon_1 -\epsilon_{n+1}$, we have $m_1(\mu +\theta) = m_1(\mu)+1$, $m_{n+1}(\mu +\theta) = m_{n+1}(\mu)-1$ and $m_j(\mu + \theta) = m_j(\mu)$ for $j \neq 1,n+1$.

If (a) is satisfied, then
\begin{align*}
\sum_{j\in J(\mu + \theta)} m_j(\mu + \theta) &= \sum_{j \in J(\mu)}m_j(\mu)+1\\
                                              &=k+1.
\end{align*}
This implies $\mu + \theta \in \wt B((k+1)\theta) \setminus \wt B(k\theta)$ by Lemma~\ref{lem:sum} (i).

If (b) is satisfied, then
\begin{align*}
\sum_{j \in J(\mu + \theta)} m_j(\mu + \theta) &= \sum_{j \in J(\mu)} m_j(\mu)+1-1\\
                                               &= \sum_{j \in J(\mu)} m_j(\mu)\\
                                               &=k.
\end{align*}
This implies $\mu + \theta \in \wt B(k\theta) \setminus \wt B((k-1)\theta)$.

If (c) is satisfied, then
\begin{align*}
\sum_{j \in J(\mu + \theta)} m_j(\mu + \theta) &= \sum_{j \in J(\mu)} m_j(\mu)\\
                                               &=k.
\end{align*}
This implies $\mu + \theta \in \wt B(k\theta) \setminus \wt B((k-1)\theta)$.

If (d) is satisfied, then
\begin{align*}
\sum_{j \in J(\mu + \theta)} m_j(\mu + \theta) &= \sum_{j \in J(\mu)} m_j(\mu)-1\\
                                               &=k-1.
\end{align*}
This implies $\mu + \theta \in \wt B((k-1)\theta) \setminus \wt B((k-2)\theta)$.
\end{proof}

Similar results for type $C_n$ and $B_n$ are verified as follows.

\begin{lem}\label{lem:cwt}
Let $\mathfrak{g}_0$ be of type $C_n$ and assume $\mu \in \wt B(k\theta) \setminus \wt B((k-1)\theta)$.
Then we have
\begin{enumerate}
\item $\mu + \theta \in \wt B((k+1)\theta) \setminus \wt B(k\theta)$ if and only if $m_1(\mu) \geq 0$.

\item $\mu + \theta \in \wt B(k\theta) \setminus \wt B((k-1)\theta)$ if and only if $m_1(\mu) = -1$.

\item $\mu + \theta \in \wt B((k-1)\theta) \setminus \wt B((k-2)\theta)$ if and only if $m_1(\mu) \leq -2$. 
\end{enumerate}					
\end{lem}

\begin{proof}
The assumption $\mu \in \wt B(k\theta) \setminus \wt B((k-1)\theta)$ implies $|\mu| = 2k$ by Lemma~\ref{lem:sum}~(ii).
Since $\theta = 2\epsilon_1$, we have 
\begin{align*}
|\mu + \theta| &= |m_1(\mu) + 2| + \sum_{j=2}^n |m_j(\mu)|\\
               &= |m_1(\mu) + 2| - |m_1(\mu)| +|\mu|\\
               &= |m_1(\mu) + 2| - |m_1(\mu)| + 2k.
\end{align*}
Then 
\begin{quote}
$|\mu + \theta| = 2(k+1)$ if and only if $m_1(\mu) \geq 0$,

$|\mu + \theta| = 2k$ if and only if $m_1(\mu) = -1$,

$|\mu + \theta| = 2(k-1)$ if and only if $m_1(\mu) \leq -2$.
\end{quote}
This completes the proof.
\end{proof}

\begin{lem}\label{lem:dwt}
Let $\mathfrak{g}_0$ be of type $B_n$ and assume $\mu \in \wt B(k\theta) \setminus \wt B((k-1)\theta)$.
Then we have
\begin{enumerate}
\item $\mu + \theta \in \wt B((k+1)\theta) \setminus \wt B(k\theta)$ if and only if $m_1(\mu) \geq 0$.
\item $\mu + \theta \in \wt B((k-1)\theta) \setminus \wt B((k-2)\theta)$ if and only if $m_1(\mu) < 0$. 
\end{enumerate}					
\end{lem}

\begin{proof}
By Lemma~\ref{lem:sum} (iii), $|\mu| = k$.
Since $\theta = \epsilon_1$, we have 
\begin{align*}
|\mu + \theta| &= |m_1(\mu) + 1| + \sum_{j=2}^n |m_j(\mu)|\\
               &= |m_1(\mu) + 1| - |m_1(\mu)| +|\mu|\\
               &= |m_1(\mu) + 1| - |m_1(\mu)| + k.
\end{align*}
Then 
\begin{quote}
$|\mu + \theta| = k+1$ if and only if $m_1(\mu) \geq 0$,

$|\mu + \theta| = k-1$ if and only if $m_1(\mu) < 0$.
\end{quote}
This completes the proof.
\end{proof}

\section{The case of type $A \sb{n} \sp{(1)}$}
We assume that $\mathfrak{g}$ is the affine Lie algebra of type $A \sb{n} \sp{(1)}$ in this section.

\subsection{The crystal structure of $B(\lambda)$}
We recall the structure of the $U_q(\mathfrak{g}_0)$-crystal $B(\lambda)$ following \cite{kn}. 
We identify the set of dominant weights $P_0^+$ with the set of Young diagrams with depth at most $n$ as follows.
For the Young diagram corresponding to a partition $(\lambda_1, \dots, \lambda_n)$, associate $\lambda = \sum_{i=1}^n (\lambda_i - \lambda_{i+1})\varpi_i \in P_0^+$ where $\lambda_{n+1}=0$.
We regard the crystal $B(\varpi_1)$ as the set consisting of letters $1,\dots,n+1$ with the following crystal structure:
\[\wt(j)=\epsilon_j,\] 
\[\e{i}(i+1)=i,\ \e{i}(j)=0\ \text{for}\ j \neq i+1,\]
\[\f{i}(i)=i+1,\ \f{i}(j)=0\ \text{for}\ j \neq i.\]
Let $\lambda \in P_0^+$.
As a set, $B(\lambda)$ is identified with the set of semistandard tableaux of shape $\lambda$ with entries $1,\dots,n+1$.
For $b\in B(\lambda)$, $\jw (b)$ denotes the Japanese reading word of $b$, that is, the word obtained by reading entries of $b$ from the top to the bottom in each column, from the right-most column to the left.
If $\jw (b) = b_1 \cdots b_N$ then we regard $b$ as $b_1 \otimes \dots \otimes b_N$, an element of $B(\varpi_1)^{\otimes N}$.
Then Kashiwara operators act on $b$ by the tensor product rule.
More explicitly, one can calculate $\e{i}b$ and $\f{i}b$ as follows.
Fix $i \in I$.
Eliminate every letter which is neither $i$ nor $i+1$ from $\jw (b)$.
If the resulting word has adjacent pairs $i \cdot (i+1)$, cancel out them, and repeat the procedure until we obtain the word $(i+1)^r \cdot i^s$ for some $r,s$.
Conclude that $\e{i} b=0$ if $r=0$ and $\f{i} b=0$ if $s=0$.
Otherwise, we determine $b_j$ and $b_k$ such that
\[\e{i}(b_1 \otimes \dots \otimes b_N) = b_1 \otimes \cdots \otimes (\e{i} b_j) \otimes \cdots \otimes b_N\]
and
\[\f{i}(b_1 \otimes \dots \otimes b_N) = b_1 \otimes \cdots \otimes (\f{i} b_k) \otimes \cdots \otimes b_N\] as $b_j$ corresponds to the right-most $i+1$ and $b_k$ the left-most $i$ in the word.   
When we deal with a tensor product $B(\lambda) \otimes B(\mu)$ for $\lambda, \mu \in P_0^+$, define $\jw (b_1 \otimes b_2) = \jw(b_1) \jw(b_2)$ for $b_1 \in B(\lambda)$ and $b_2 \in B(\mu)$.
Then one can calculate the actions of Kashiwara operators in the same way. 

\subsection{Kirillov-Reshetikhin crystals for type $A$}

We recall the crystal structure of $B^{i,l}$ following \cite{shi}.
Since $W^{i,l}$ is isomorphic to $V(l\varpi_i)$ as a $U_q(\mathfrak{g}_0)$-module, $B^{i,l}$ is isomorphic to $B(l\varpi_i)$ as a $U_q(\mathfrak{g}_0)$-crystal.
For the description of the actions of $\e{0}$ and $\f{0}$, we use the promotion operator $\sigma$.
The operator $\sigma$ is a bijection from $B^{i,l}$ to $B^{i,l}$ which satisfy
\[\sigma \e{j} = \e{j+1}\sigma\]
and   
\[\sigma \f{j} = \f{j+1}\sigma\]
for each $j$ modulo $n+1$.
This map corresponds to the automorphism of the Dynkin diagram of type $A \sb{n} \sp{(1)}$ which takes $j$ to $j+1$ modulo $n+1$.
Using $\sigma$, we obtain
\[\e{0} = \sigma^{-1}\e{1}\sigma\]
and   
\[\f{0} = \sigma^{-1}\f{1}\sigma.\]
We shall describe the action of $\sigma$ explicitly only for $B^{1,l}$ and $B^{n,l}$, and omit the general definition.
See \cite[3.3]{shi} for details.  

For $b \in B^{1,l}$, we denote the number of entries $j$ appearing in the tableau $b$ by $x_j(b)$.
The map which takes $b \in B^{1,l}$ to $(x_j(b))_{j = 1, \dots, n+1} \in (\mathbb{Z}_{\geq 0}) \sp{n+1}$ is injective.
The action of the promotion operator is given by
\[x_1(\sigma(b)) = x_{n+1}(b)\]
and
\[x_j(\sigma(b)) = x_{j-1}(b)\]
for $j=2,\dots,n+1$.
Hence the actions of Kashiwara operators are described as follows:
\begin{align*}
x_j(\e{0} b) & = x_j(b) - \delta_{j1} + \delta_{j \, n+1},\\
x_j(\f{0} b) & = x_j(b) + \delta_{j1} - \delta_{j \, n+1},\\
x_j(\e{i} b) &= x_j(b) - \delta_{j \, i+1} + \delta_{ji} \quad \text{for}\ i \in I_0,\\
x_j(\f{i} b) &= x_j(b) + \delta_{j \, i+1} - \delta_{ji} \quad \text{for}\ i \in I_0.
\end{align*}
Here we define $b=0$ if $x_j(b)<0$ for some $j$.
Immediately we have
\begin{align*}
\varepsilon_0(b) & = x_1(b),\\
\varphi_0(b) & = x_{n+1}(b),\\
\varepsilon_i(b) & = x_{i+1}(b) \quad \text{for}\ i \in I_0,\\
\varphi_i(b) & = x_{i}(b) \quad \text{for}\ i \in I_0.
\end{align*}

For $j= 1,\dots,n+1$, we denote by $C_j$ the semistandard tableau consisting of one column with depth $n$ which has no entry $j$.
Each column in a semistandard tableau of shape $(l^n)$ is $C_j$ for some $j$.
These columns are arranged as, from left to right, some (maybe 0) $C_{n+1}$'s, some $C_{n}$'s, $\dots$, and some $C_{1}$'s. 
For $b \in B^{n,l}$, we denote the number of columns $C_j$ in $b$ by $y_j(b)$.
The map which takes $b \in B^{n,l}$ to $(y_j(b))_{j = 1, \dots, n+1} \in (\mathbb{Z} \sb{\geq 0}) \sp{n+1}$ is injective.
The action of the promotion operator is given by
\[y_1(\sigma(b)) = y_{n+1}(b)\]
and
\[y_j(\sigma(b)) = y_{j-1}(b)\]
for $j=2,\dots,n+1$.
Hence the actions of Kashiwara operators are described as follows:
\begin{align*}
y_j(\e{0}b) &= y_j(b) - \delta_{j\,n+1} + \delta_{j1},\\
y_j(\f{0}b) &= y_j(b) + \delta_{j\,n+1} - \delta_{j1},\\
y_j(\e{i}b) &= y_j(b) - \delta_{ji} + \delta_{j\,i+1} \quad \text{for}\ i \in I_0,\\
y_j(\f{i}b) &= y_j(b) + \delta_{ji} - \delta_{j\,i+1} \quad \text{for}\ i \in I_0.
\end{align*}
Here we define $b=0$ if $y_j(b)<0$ for some $j$.
Immediately we have
\begin{align*}
\varepsilon_0(b) &= y_{n+1}(b),\\
\varphi_0(b) &= y_1(b),\\
\varepsilon_i(b) &= y_{i}(b) \quad \text{for}\ i \in I_0,\\
\varphi_i(b) &= y_{i+1}(b) \quad \text{for}\ i \in I_0.
\end{align*}

\subsection{The crystal structure of $B_l$}

For a fixed nonnegative integer $l$, we set $B_l=B^{1,l}\otimes B^{n,l}$.
We shall describe the decomposition of the tensor product $B(l\varpi_1) \otimes B(l\varpi_n)$.
Let $b_1$ be a semistandard tableau of shape $(k)$ and $b_2$ that of shape $(k^n)$, and assume $x_1(b_1) = 0$ or $y_1(b_2) = 0$.
We define $b_1 \cdot b_2$ as the semistandard tableau of shape $(2k,k^{n-1})$ such that $\jw (b_1 \cdot b_2) = \jw (b_1) \jw (b_2)$. 
Define a map $\alpha \colon B(l\varpi_1) \otimes B(l\varpi_n) \to \bigoplus_{k=0}^l B(k\theta)$ by $\alpha (b_1 \otimes b_2) = \tilde{b}_1 \cdot \tilde{b}_2 \in B(k\theta)$ where $k = l - \min \{x_1(b_1),y_1(b_2)\}$, $\tilde{b}_1$ is the semistandard tableau of shape $(k)$ which is obtained by removing $\min \{x_1(b_1),y_1(b_2)\}$ $1$'s from $b_1$, and $\tilde{b}_2$ is the semistandard tableau of shape $(k^n)$ obtained by removing $\min \{x_1(b_1),y_1(b_2)\}$ $C_1$'s from $b_2$. 
Then we have the following lemma. 

\begin{lem}\label{lem:isom}
The map $\alpha$ is an isomorphism of $U_q(\mathfrak{g}_0)$-crystals between $B(l\varpi_1) \otimes B(l\varpi_n)$ and $\bigoplus_{k=0}^l B(k\theta)$.
\end{lem}

\begin{proof}
It is obvious that $\alpha$ is bijective and preserves weights.
Suppose $\alpha (b_1 \otimes b_2) = \tilde{b}_1 \cdot \tilde{b}_2 \in B(k\theta)$.
Then we have $\jw (b_1 \otimes b_2) = \jw(\tilde{b}_1) \cdot 1^{l-k} \cdot (2 \cdot 3 \cdots (n+1))^{l-k} \cdot \jw(\tilde{b}_2)$.
Since the element $1 \otimes 2 \otimes \dots \otimes (n+1)$ is annihilated by all $\e{i}$ and $\f{i}$ for $i \in I_0$, 
$\alpha$ commutes with Kashiwara operators. 
\end{proof}

We write $b \in B(k\theta) \subset B_l$ for $b \in B_l$, if $\alpha(b) \in B(k\theta)$.
Here $B_l$ is identified with $B(l\varpi_1) \otimes B(l\varpi_n)$ as a set.

For $j=1,\dots,n+1$, we define an injective map $\Theta_j\colon B_{l-1}\to B_l$ as follows.
For $b_1\otimes b_2\in B_l$, set $\Theta_j(b_1\otimes b_2)=b_1'\otimes b_2'$, where $b_1'$ is the semistandard tableau of shape ($l$) which is obtained by adding an entry $j$ to $b_1$, and $b_2'$ is the semistandard tableau of shape ($l^n$) obtained by adding a column $C_j$ to $b_2$. 
These maps preserve $U_q(\mathfrak{g}_0)$-weights by definition.

\begin{prop}\label{prop:1}
\begin{enumerate}
\item For $k=0,\dots,l-1$, $\Theta_1(B(k\theta))\subset B(k\theta)$.
\item For $i \in I_0$, the map $\Theta_1$ and the Kashiwara operator $\f{i}$ commute with each other.
\item Let $b\in B_{l-1}$ and assume $\f{0} b\neq 0$. Then $\Theta_1(\f{0} b)=\f{0}\Theta_1(b)$. 
\item Let $b\in B_{l-1}$ and assume $\f{0} b=0$. Then $\varphi_0(\Theta_1(b))=1$ and $\f{0}\Theta_1(b)\in B(l\theta)\subset B_l$.
\end{enumerate}
\end{prop}

\begin{proof}
Let $b_1 \otimes b_2 \in B_{l-1}$ and suppose $\Theta_1(b_1\otimes b_2)=b_1'\otimes b_2'$.
We have $x_1(b_1')=x_1(b_1)+1$ and $y_1(b_2')=y_1(b_2)+1$ by the definition of $\Theta_1$.
This implies (i).

Since $\jw (b_1' \otimes b_2') = \jw (b_1) \cdot 1 \cdot 2 \cdots (n+1) \cdot \jw (b_2)$, (ii) is immediate.

We prove (iii).
By the definition of $\Theta_1$, we have 
\[\varphi_0(b_1') = x_{n+1}(b_1') = x_{n+1}(b_1) = \varphi_0(b_1)\]
and
\[\varepsilon_0(b_2') = y_{n+1}(b_2') = y_{n+1}(b_2) =  \varepsilon_0(b_2).\]
There are the following two cases:
\[\f{0} (b_1 \otimes b_2) =
\begin{cases}
(\f{0}b_1) \otimes b_2 & \text{if} \ \varphi_0 (b_1) > \varepsilon_0 (b_2),\\
b_1 \otimes (\f{0}b_2) & \text{if} \ \varphi_0 (b_1) \leq \varepsilon_0 (b_2).
\end{cases}\]
First we assume $\f{0} (b_1 \otimes b_2) = (\f{0} b_1) \otimes b_2$.
Then $\f{0} (b_1' \otimes b_2') = (\f{0} b_1') \otimes b_2'$ since $\varphi_0(b_1') = \varphi_0(b_1)$ and $\varepsilon_0(b_2') = \varepsilon_0(b_2)$. 
Recalling the action of $\f{0}$, the tableau $\f{0} b_1$ (resp.\ $\f{0} b_1'$) is obtained from $b_1$ (resp.\ $b_1'$) by removing an entry $n+1$ and adding an entry $1$.
We have
\begin{align*}
\Theta_1 (\f{0} (b_1 \otimes b_2)) &= \Theta_1 ((\f{0} b_1) \otimes b_2) \\
                                   &= (\f{0} b_1)' \otimes b_2'
\end{align*}
where $(\f{0} b_1)'$ is the tableau obtained from $\f{0} b_1$ by adding an entry $1$.
Hence the assertion is true since $(\f{0} b_1)'$ coincides with $\f{0} b_1'$: the semistandard tableau of shape ($l$) which is obtained from $b_1$ by removing one entry $n+1$ and adding two $1$'s.  
Next we assume $\f{0} (b_1 \otimes b_2) = b_1 \otimes (\f{0} b_2)$.
Similarly, we have $\f{0} (b_1' \otimes b_2') = b_1' \otimes (\f{0} b_2')$.
In this case, we have
\begin{align*}
\Theta_1 (\f{0} (b_1 \otimes b_2)) &= \Theta_1 ( b_1 \otimes (\f{0} b_2)) \\
                                   &= b_1' \otimes (\f{0} b_2)'
\end{align*}
where $(\f{0} b_2)'$, unless $\f{0} b_2 = 0$, is obtained from $b_2$ by removing one column $C_1$, adding one $C_{n+1}$, and adding one $C_1$ in this order.
It coincides with $\f{0} b_2'$ whenever $b_2$ has at least one column $C_1$.
Thus (iii) is proved.

To prove (iv), suppose $b = b_1 \otimes b_2 \in B_{l-1}$ and $\f{0}b=0$.
We have $\varphi_0(b_2)=0$ and $\varphi_0(b_1)\leq \varepsilon_0(b_2)$ by the tensor product rule.
Then $\f{0}(b_1'\otimes b_2')=b_1'\otimes (\f{0}b_2')$ and $\f{0}b_2'\neq 0$ since $b_2'$ has exactly one column $C_1$.
As $\f{0}b_2'$ has no more $C_1$, $b_1'\otimes (\f{0}b_2')\in B(l\theta)\subset B_l$ and $\f{0}^2(b_1'\otimes b_2')=0$.
\end{proof}

By Proposition~\ref{prop:1}, the crystal graph of $B_{l-1}$ can be naturally regarded as a full subgraph of that of $B_l$.

\begin{prop}\label{prop:2}
Let $j = 2,\dots,n+1$.
\begin{enumerate}
\item For $k = 0,\dots,l-1$, $\Theta_j(B(k\theta)) \subset B((k+1)\theta)$.
\item Let $i \in I_0$ and $b \in B_{l-1}$.
If $\f{i} b \neq 0$, then $\Theta_j(\f{i} b) = \f{i}\Theta_j(b)$.
\item The map $\Theta_j$ and the Kashiwara operator $\f{0}$ commute with each other.
\end{enumerate}
\end{prop}

\begin{proof}
For $b \in B_{l-1}$, we suppose $b = b_1 \otimes b_2$ and $\Theta_j(b) = b_1' \otimes b_2'$.
Then (i) follows from $x_1(b_1') = x_1(b_1)$ and $y_1(b_2') = y_1(b_2)$.

Recalling $\varphi_i (b_1) = x_i(b_1)$, $\varepsilon_i(b_2) = y_i(b_2)$ etc., we have
\[\varphi_i (b_1') = \varphi_i (b_1) + \delta_{ij}\ \text{for}\ i \in I_0,\]
\[\varepsilon_i (b_2') = \varepsilon_i (b_2) + \delta_{ij}\ \text{for}\	 i \in I_0,\]
\[\varphi_0 (b_1') = \varphi_0 (b_1) + \delta_{j\,n+1},\]
\[\varepsilon_0 (b_2') = \varepsilon_0 (b_2) + \delta_{j\, n+1}.\] 
Then (ii) and the case $\f{0}b \neq 0$ of (iii) are obtained from the above formulas by arguments similar to that in the proof of Proposition \ref{prop:1} (iii).
In the case $\f{0}b = 0$, we have $\f{0}(b_1'\otimes b_2') = b_1'\otimes (\f{0}b_2') = 0$ since $\varphi_0(b_2') = y_1(b_2') = y_1(b_2) = \varphi_0(b_2)=0$. 
\end{proof}

Hence the action of $\f{0}$ on the set $\bigcup_{j=2}^{n+1} \Ima \Theta_j \subset B_l$ is completely determined if we know that on $B_{l-1}$.

\begin{rem}\label{rem}
For Proposition~\ref{prop:1} and \ref{prop:2}, the similar statements hold when we replace $\f{i}$ by $\e{i}$ for $i \in I$ and $\varphi_0$ by $\varepsilon_0$.
We can prove them in a similar way.  
\end{rem}

From now on, we study the action of $\f{0}$ on $B_l \setminus \bigcup_{j=2}^{n+1} \Ima \Theta_j$.
Let $\mu \in P_0$. 
In Subsection~\ref{subsec:lemma}, we defined $m_j(\mu)$ as the coefficient of $\epsilon_j$ in $\mu$ with $\sum_{j=1}^{n+1} m_j(\mu) = 0$. 
The number of entries $j$ in a semistandard tableau of shape $k\theta$ and weight $\mu$ is equal to $m_j(\mu) + k$.
Put $m_j(\mu, k) = m_j(\mu) + k$.
Then we can rewrite Lemma~\ref{lem:sum} (i) as follows.

\begin{lem}\label{lem:sum2}
Let $k$ be a nonnegative integer.
For $\mu \in P_0$, $\mu \in \wt B(k\theta) \setminus \wt B((k-1)\theta)$ if and only if $\sum_{j \in J(\mu)}  m_j(\mu, k) = k|J(\mu)| + k$.
\end{lem}

\begin{prop}\label{prop:set} 
We have
\[\bigcup_{j=2}^{n+1} \Theta_j(B((k-1)\theta)) = \{b \in B(k\theta) \mid \wt b \in \wt B((k-1)\theta)\}.\]
\end{prop}

\begin{proof}
Proposition~\ref{prop:2} (i) implies that the right-hand side contains the left.
Suppose $b$ does not belong to the left-hand side and $\alpha(b) = \tilde{b}_1 \cdot \tilde{b}_2 \in B(k\theta)$.
Put $\mu = \wt b$.
Let $J$ be the set of the letters appearing in $\tilde{b}_1$. 
Then $\tilde{b}_2$ has no column $C_j$ if $j \in J$.
This means the number of entries $j$ in $\tilde{b}_2$ is equal to $k$ if $j \in J$.
Therefore $J(\mu) = J$ and $\sum_{j \in J} m_j(\mu ,k) = k|J|+k$.
We obtain $\mu \in \wt B(k\theta) \setminus \wt B((k-1)\theta)$ by Lemma~\ref{lem:sum2}.
\end{proof}

Although the following fact seems to be known, we give a proof. 

\begin{lem}\label{lem:free}
Let $k$ be a nonnegative integer.
The multiplicity of every element of $\wt B(k\theta) \setminus \wt B((k-1)\theta)$ in $B(k\theta)$ is one. 
\end{lem}

\begin{proof}
Suppose $b = b_1 \cdot b_2 \in B(k\theta)$ where $b_1 \in B(k\varpi_1)$ and $b_2 \in B(k\varpi_n)$, and $\mu = \wt b \in \wt B(k\theta) \setminus \wt B((k-1)\theta)$.
By Lemma \ref{lem:sum2}, the number of boxes of $b$ filled with elements of $J(\mu)$ is equal to $k|J(\mu)|+k$.
Therefore $b_1$ should have only elements of $J(\mu)$ with its entries and the number of entries $j$ in $b_2$ should be equal to $k$ for $j \in J(\mu)$.
When the weight of $b$ is given, the numbering of $b$ is uniquely determined by that of $b_1$.
Hence $b$ is uniquely determined by its weight.
\end{proof}

The following proposition immediately follows from Proposition~\ref{prop:set} and Lemma~\ref{lem:free}.

\begin{prop}\label{prop:rest}
The restriction $\wt \mid_{B_l \setminus \bigcup_{j=2}^{n+1} \Ima \Theta_j}$ is injective.
\end{prop}

Let $b \in B(k\theta) \subset B_l$ and $b \notin \bigcup_{j=2}^{n+1}\Ima \Theta_j$.
By Lemma~\ref{lem:wt}, exactly one of the following occurs:
\begin{itemize}
\item $\wt b + \theta \in \wt B((k+1)\theta) \setminus \wt B(k\theta)$.

\item $\wt b + \theta \in \wt B(k\theta) \setminus \wt B((k-1)\theta)$.

\item $\wt b + \theta \in \wt B((k-1)\theta) \setminus \wt B((k-2)\theta)$.
\end{itemize}

\begin{thm}\label{thm:main}
Let $b \in B(k\theta) \subset B_l$ and assume $b \notin \bigcup_{j=2}^{n+1}\Ima \Theta_j$.
\begin{enumerate}
\item We have $\f{0} b=0$ if and only if $k=l$ and $\wt b + \theta \notin \wt B(l\theta)$. 

\item If $\wt b + \theta \in \wt B((k+1)\theta) \setminus \wt B(k\theta)$, then $\f{0} b \in B((k+1)\theta) \subset B_l$.

\item If $\wt b + \theta \in \wt B(k\theta) \setminus \wt B((k-1)\theta)$, then $\f{0} b \in B(k\theta) \subset B_l$.

\item If $\wt b + \theta \in \wt B((k-1)\theta) \setminus \wt B((k-2)\theta)$, then $\f{0} b \in B((k-1)\theta) \subset B_l$.
\end{enumerate}
Moreover, the element $\f{0}b$ is uniquely determined by its weight in each case.
\end{thm}

\begin{proof}
We prove (i).
Since $\alpha_0 + \theta =0$ in the weight lattice $P_{\text{cl}}$, $\f{0}$ raises $U_q(\mathfrak{g}_0)$-weights by $\theta$.
Hence the sufficiency is obvious.
Suppose $b = b_1 \otimes b_2$.
If $\f{0} b = 0$, we have $\varphi_0(b_2) = 0$ and $\varphi_0(b_1) \leq \varepsilon_0(b_2)$ by the tensor product rule.
Since $y_1(b_2) = \varphi_0(b_2) = 0$, we have $k=l$.
The number of entries $1$ in $b$ is 
\begin{align*}
x_1(b_1) + (l - y_1(b_2)) &= l + x_1(b_1)\\
                          &\geq l
\end{align*}
and the number of entries $n+1$ is
\begin{align*}
x_{n+1}(b_1) + (l - y_{n+1}(b_2)) &= l + (\varphi_0(b_1) - \varepsilon_0(b_2)) \\
                                  & \leq l.
\end{align*}
Hence we have $m_1(\wt b) \geq 0$ and $m_{n+1}(\wt b) \leq 0$.
By Lemma~\ref{lem:wt}~(i), we obtain $\wt b + \theta \notin \wt B(l\theta)$.

By (i), $\wt b + \theta \in \wt B(l\theta)$ implies $\f{0}b \neq 0$.
Then we see $\f{0}b \in B_l \setminus \bigcup_{j=2}^{n+1} \Ima \Theta_j$ by Proposition~\ref{prop:2}~(iii) and Remark~\ref{rem}.
By Proposition~\ref{prop:rest}, $\f{0}b$ is uniquely determined by its weight $\wt b + \theta$.
This proves the remaining assertions.
\end{proof}

\section{The case of $C \sb{n} \sp{(1)}$}

Assume that $\mathfrak{g}$ is of type $C \sb{n} \sp{(1)}$ in this section.
The crystal structure of $B_l = B^{1,2l}$ is given in \cite{kkm}.
We recall explicit formulas on the actions of Kashiwara operators in \cite{okadomem}.
As sets,
\[B(k\theta) = \{(x_1, \dots, x_n, \bar{x}_n, \dots, \bar{x}_1) \in (\mathbb{Z} \sb{\geq 0})^{2n} \mid \sum_{j=1}^n (x_j + \bar{x}_j) = 2k\}\]
and
\[B_l = \bigsqcup_{k=0}^l B(k\theta).\] 
For $b = (x_1, \dots, \bar{x}_1) \in B_l$, the $U_q(\mathfrak{g}_0)$-weight of $b$ and the actions of Kashiwara operators are given as follows:
\[\wt b = \sum_{j=1}^n (x_j - \bar{x}_j)\epsilon_j,\]
\begin{align*}
\e{0} b &=
\begin{cases}
(x_1 - 2, x_2, \dots, \bar{x}_2, \bar{x}_1) & \text{ if } x_1 \geq \bar{x}_1 + 2,\\
(x_1 - 1, x_2, \dots, \bar{x}_2, \bar{x}_1 + 1) & \text{ if } x_1 = \bar{x}_1 + 1,\\
(x_1, x_2, \dots, \bar{x}_2, \bar{x}_1 + 2) & \text{ if } x_1 \leq \bar{x}_1,
\end{cases}\\
\e{i} b &=
\begin{cases}
(x_1, \dots, x_i + 1, x_{i+1} - 1, \dots, \bar{x}_1) & \text{ if } x_{i+1} > \bar{x}_{i+1},\\
(x_1, \dots, \bar{x}_{i+1} + 1, \bar{x}_i - 1, \dots, \bar{x}_1) & \text{ if } x_{i+1} \leq \bar{x}_{i+1},\\
\end{cases}\\
\e{n} b &= (x_1, \dots, x_n + 1, \bar{x}_n - 1, \dots, \bar{x}_1),\\
\f{0} b &=
\begin{cases}
(x_1 + 2, x_2, \dots, \bar{x}_2, \bar{x}_1) & \text{ if } x_1 \geq \bar{x}_1,\\
(x_1 + 1, x_2, \dots, \bar{x}_2, \bar{x}_1 - 1) & \text{ if } x_1 = \bar{x}_1 - 1,\\
(x_1, x_2, \dots, \bar{x}_2, \bar{x}_1 - 2) & \text{ if } x_1 \leq \bar{x}_1 - 2,
\end{cases}\\
\f{i} b &=
\begin{cases}
(x_1, \dots, x_i - 1, x_{i+1} + 1, \dots, \bar{x}_1) & \text{ if } x_{i+1} \geq \bar{x}_{i+1},\\
(x_1, \dots, \bar{x}_{i+1} - 1, \bar{x}_i + 1, \dots, \bar{x}_1) & \text{ if } x_{i+1} < \bar{x}_{i+1},\\
\end{cases}\\
\f{n} b &= (x_1, \dots, x_n - 1, \bar{x}_n + 1, \dots, \bar{x}_1).
\end{align*}
In each case, if the right-hand side does not belong to the set $B_l$, regard it as zero.
It is obvious that the crystal graph of $B_{l-1}$ is a full subgraph of that of $B_l$ from the above formulas.
For $b \in B(k\theta) \subset B_l$, we have
\begin{align*}
\varepsilon_0(b) &= (l-k) + \max(0, x_1 - \bar{x}_1),\\
\varphi_0(b) &= (l-k) + \max(0, \bar{x}_1 - x_1),\\
\varepsilon_i(b) &= \bar{x}_i + \max(0, x_{i+1} - \bar{x}_{i+1}),\\
\varphi_i(b) &= x_i + \max(0, \bar{x}_{i+1} - x_{i+1}),\\
\varepsilon_n(b) &= \bar{x}_n,\\
\varphi_n(b) &= x_n.
\end{align*}

For $j = 1, \dots, n$, we define a map $\Phi_j \colon B_{l-1} \to B_l$ by
\[\Phi_j(x_1, \dots, \bar{x}_1) = (\dots, x_j + 1, \dots, \bar{x}_j +1, \dots).\] 
\begin{prop}\label{prop:commutec}
Let $j=1, \dots, n$.
\begin{enumerate}
\item For $k = 0, \dots, l-1$, $\Phi_j(B(k\theta)) \subset B((k+1)\theta)$.

\item Let $i \in I_0$ and $b \in B_{l-1}$.
If $\f{i} b \neq 0$, then $\Phi_j(\f{i} b) = \f{i} \Phi_j(b)$.

\item The map $\Phi_j$ and the Kashiwara operator $\f{0}$ commute with each other.
\end{enumerate}
\end{prop}

\begin{proof}
The assertion of (i) is obvious from the definition of $\Phi_j$.

Suppose $b = (x_1, \dots, \bar{x}_1)$ and $\Phi_j (b) = (y_1, \dots, \bar{y}_1)$.
Then $y_i - \bar{y}_i = x_i - \bar{x}_i$ for $i = 1, \dots, n$.
Hence (ii) and the case $\f{0} b \neq 0$ of (iii) are immediate from the formulas on the actions of Kashiwara operators.
If $\f{0} b = 0$, then $b \in B((l-1)\theta) \subset B_{l-1}$ and $x_1 \geq \bar{x}_1$ by the formula on $\varphi_0$.
Hence $\Phi_j(b) \in B(l\theta)$ and $\f{0} \Phi_j(b) = 0$.
\end{proof}

\begin{prop}\label{prop:setc} 
We have
\[\bigcup_{j=1}^n \Phi_j(B((k-1)\theta)) = \{b \in B(k\theta) \mid \wt b \in \wt B((k-1)\theta)\}.\]
\end{prop}

\begin{proof}
Let $b = (x_1, \dots, \bar{x}_1) \in B(k\theta)$.
By Lemma~\ref{lem:sum}~(ii), $\wt b \in \wt B(k\theta) \setminus \wt B((k-1)\theta)$ if and only if $x_j = 0$ or $\bar{x}_j = 0$
for each $j = 1, \dots, n$.
This condition is equivalent to that $b$ does not belong to the left-hand side.
\end{proof}

\begin{lem}\label{lem:freec}
Let $k$ be a nonnegative integer.
The multiplicity of every element of $\wt B(k\theta) \setminus \wt B((k-1)\theta)$ in $B(k\theta)$ is one. 
\end{lem}

\begin{proof}
Let $b = (x_1, \dots, \bar{x}_1) \in B(k\theta)$ and $\mu = \wt b \in \wt B(k\theta) \setminus \wt B((k-1)\theta)$.
Set
\[J_+ = \{ j \mid m_j(\mu) > 0\},\]
\[J_- = \{ j \mid m_j(\mu) < 0\}.\]
Then we have 
\begin{align*}
x_j  &=
\begin{cases}
m_j(\mu) & \text{ for } j \in J_+,\\
0   & \text{ for } j \notin J_+,
\end{cases}\\
\bar{x}_j &=
\begin{cases}
-m_j(\mu) & \text{ for } j \in J_-,\\
0   & \text{ for } j \notin J_-
\end{cases}
\end{align*}
by Proposition~\ref{prop:setc}.
This means that $b$ is uniquely determined by its weight.
\end{proof}

\begin{thm}\label{thm:mainc}
Let $b \in B(k\theta) \subset B_l$ and assume $b \notin \bigcup_{j=1}^n \Ima \Phi_j$.
\begin{enumerate}
\item We have $\f{0} b=0$ if and only if $k=l$ and $\wt b + \theta \notin \wt B(l\theta)$. 

\item If $\wt b + \theta \in \wt B((k+1)\theta) \setminus \wt B(k\theta)$, then $\f{0} b \in B((k+1)\theta) \subset B_l$.

\item If $\wt b + \theta \in \wt B(k\theta) \setminus \wt B((k-1)\theta)$, then $\f{0} b \in B(k\theta) \subset B_l$.

\item If $\wt b + \theta \in \wt B((k-1)\theta) \setminus \wt B((k-2)\theta)$, then $\f{0} b \in B((k-1)\theta) \subset B_l$.
\end{enumerate}
Moreover, the element $\f{0}b$ is uniquely determined by its weight in each case.
\end{thm}
\begin{proof}
Suppose $b = (x_1, \dots, \bar{x}_1)$.
Then $m_1(\wt b) = x_1 - \bar{x}_1$.
Hence the assertions are immediate from the formula on the action of $\f{0}$, Lemma~\ref{lem:cwt} and Lemma~\ref{lem:freec}. 
\end{proof}

\section{The case of $D \sb{n+1} \sp{(2)}$}

Assume that $\mathfrak{g}$ is of type $D \sb{n+1} \sp{(2)}$ in this section.
Thus $\mathfrak{g}_0$ is of type $B_n$.
The crystal structure of $B_l = B^{1,l}$ is given in \cite{kmn}.
Explicit formulas are available in \cite{okadomem} and we follow it.
As sets,
\begin{align*}
B(k\theta) = \{(x_1, \dots, x_n, x_0, \bar{x}_n, \dots, \bar{x}_1) & \in (\mathbb{Z} \sb{\geq 0})^{2n + 1} \mid \\
           & x_0 = 0 \text{ or } 1, \sum_{j=1}^n (x_j + \bar{x}_j) + x_0 = k\}
\end{align*}
and
\[B_l = \bigsqcup_{k=0}^l B(k\theta).\] 
For $b = (x_1, \dots, \bar{x}_1) \in B_l$, the $U_q(\mathfrak{g}_0)$-weight of $b$ and the actions of Kashiwara operators are given as follows:
\[\wt b = \sum_{j=1}^n (x_j - \bar{x}_j)\epsilon_j,\]
\begin{align*}
\e{0} b &=
\begin{cases}
(x_1 - 1, x_2, \dots, \bar{x}_2, \bar{x}_1) & \text{ if } x_1 > \bar{x}_1,\\
(x_1, x_2, \dots, \bar{x}_2, \bar{x}_1 + 1) & \text{ if } x_1 \leq \bar{x}_1,
\end{cases}\\
\e{i} b &=
\begin{cases}
(x_1, \dots, x_i + 1, x_{i+1} - 1, \dots, \bar{x}_1) & \text{ if } x_{i+1} > \bar{x}_{i+1},\\
(x_1, \dots, \bar{x}_{i+1} + 1, \bar{x}_i - 1, \dots, \bar{x}_1) & \text{ if } x_{i+1} \leq \bar{x}_{i+1},\\
\end{cases}\\
\e{n} b &=
\begin{cases}
(x_1, \dots, x_n, x_0 + 1, \bar{x}_n - 1, \dots, \bar{x}_1) & \text{ if } x_0 = 0,\\
(x_1, \dots, x_n + 1, x_0 - 1, \bar{x}_n, \dots, \bar{x}_1) & \text{ if } x_0 = 1,
\end{cases}\\
\f{0} b &=
\begin{cases}
(x_1 + 1, x_2, \dots, \bar{x}_2, \bar{x}_1) & \text{ if } x_1 \geq \bar{x}_1,\\
(x_1, x_2, \dots, \bar{x}_2, \bar{x}_1 - 1) & \text{ if } x_1 < \bar{x}_1,
\end{cases}\\
\f{i} b &=
\begin{cases}
(x_1, \dots, x_i - 1, x_{i+1} + 1, \dots, \bar{x}_1) & \text{ if } x_{i+1} \geq \bar{x}_{i+1},\\
(x_1, \dots, \bar{x}_{i+1} - 1, \bar{x}_i + 1, \dots, \bar{x}_1) & \text{ if } x_{i+1} < \bar{x}_{i+1},\\
\end{cases}\\
\f{n} b &=
\begin{cases}
(x_1, \dots, x_n - 1, x_0 + 1, \bar{x}_n, \dots, \bar{x}_1) & \text{ if } x_0 = 0,\\ 
(x_1, \dots, x_n, x_0 - 1, \bar{x}_n + 1, \dots, \bar{x}_1) & \text{ if } x_0 = 1.
\end{cases}
\end{align*}
In each case, if the right-hand side does not belong to $B_l$, regard it as zero.
The crystal graph of $B_{l-1}$ is a full subgraph of that of $B_l$.
For $b \in B(k\theta) \subset B_l$, we have
\begin{align*}
\varepsilon_0(b) &= (l-k) + 2\max(0, x_1 - \bar{x}_1),\\
\varphi_0(b) &= (l-k) + 2\max(0, \bar{x}_1 - x_1),\\
\varepsilon_i(b) &= \bar{x}_i + \max(0, x_{i+1} - \bar{x}_{i+1}),\\
\varphi_i(b) &= x_i + \max(0, \bar{x}_{i+1} - x_{i+1}),\\
\varepsilon_n(b) &= 2\bar{x}_n + x_0,\\
\varphi_n(b) &= 2x_n + x_0.
\end{align*}

We define maps $\Psi_j \colon B_{l-2} \to B_l$ for $j=1, \dots, n-1$ and $\Psi_n \colon B_{l-1} \to B_l$ by
\begin{align*}
\Psi_j(x_1, \dots, \bar{x}_1) &= (\dots, x_j + 1, \dots, \bar{x}_j + 1, \dots) \ \text{ for } j = 1, \dots, n-1,\\
\Psi_n(x_1, \dots, \bar{x}_1) &=
\begin{cases}
(x_1, \dots, x_n, x_0+1, \bar{x}_n, \dots, \bar{x}_1) \ \text{ if } x_0 = 0,\\
(x_1, \dots, x_n + 1, x_0 - 1, \bar{x}_n + 1, \dots, \bar{x}_1) \ \text{ if } x_0 = 1.
\end{cases}
\end{align*}

\begin{prop}\label{prop:commuted}
\begin{enumerate}
\item Let $j =1, \dots, n-1$.
Then $\Psi_j(B(k\theta)) \subset B((k+2)\theta)$ for $k = 0, \dots, l-2$.

\item For $k = 0, \dots, l-1$, $\Psi_n(B(k\theta)) \subset B((k+1)\theta)$.

\item Let $j = 1, \dots, n-1$ $(\text{resp.\ }j = n)$ and $i \in I_0$.
If we take $b \in B_{l-2}$ $(\text{resp.\ }b \in B_{l-1})$ with $\f{i} b \neq 0$, then $\Psi_j(\f{i} b) = \f{i} \Psi_j(b)$.

\item For $j = 1, \dots, n$, the map $\Psi_j$ and the Kashiwara operator $\f{0}$ commute with each other.
\end{enumerate}
\end{prop}

\begin{proof}
The assertions of (i) and (ii) are obvious.

The assertions of (iii) and (iv) are immediate except for the commutativity of $\Psi_n$ and $\f{n}$.
We show that $\Psi_n(\f{n} b) = \f{n} \Psi_n(b)$ when $\f{n} b \neq 0$.
Let $b = (x_1, \dots, \bar{x}_1)$.
First we assume $x_0 = 0$.
Then
\begin{align*}
\Psi_n(\f{n} b) &= \Psi_n(x_1, \dots, x_n - 1,x_0 + 1, \bar{x}_n, \dots, \bar{x}_1)\\
                &=(x_1, \dots, x_n, x_0, \bar{x}_n + 1, \dots, \bar{x}_1).
\end{align*} 
On the other hand,
\begin{align*}
\f{n}\Psi_n(b) &= \f{n}(x_1, \dots, x_n ,x_0 + 1, \bar{x}_n, \dots, \bar{x}_1)\\
                &=(x_1, \dots, x_n, x_0, \bar{x}_n + 1, \dots, \bar{x}_1).
\end{align*} 
Hence the assertion is true.
Next assume $x_0 = 1$.
Then we obtain
\begin{align*}
\Psi_n(\f{n} b) &= \Psi_n(x_1, \dots, x_n, x_0 - 1, \bar{x}_n + 1, \dots, \bar{x}_1)\\
                &=(x_1, \dots, x_n, x_0, \bar{x}_n + 1, \dots, \bar{x}_1)
\end{align*} 
and
\begin{align*}
\f{n}\Psi_n(b) &= \f{n}(x_1, \dots, x_n + 1,x_0 - 1, \bar{x}_n + 1, \dots, \bar{x}_1)\\
               &= (x_1, \dots, x_n, x_0, \bar{x}_n + 1, \dots, \bar{x}_1).
\end{align*}
\end{proof}

\begin{prop}\label{prop:dset} 
We have
\begin{align*}
\left( \bigcup_{j=1}^{n-1} \Psi_j(B((k-2)\theta)) \right) \cup\, & \Psi_n (B((k-1)\theta)) \\
& = \{ b \in B(k\theta) \mid \wt b \in \wt B((k-1)\theta) \}.
\end{align*}
\end{prop}

\begin{proof}
Let $b = (x_1, \dots, \bar{x}_1) \in B(k\theta)$.
By Lemma~\ref{lem:sum}~(iii), $\wt b \in \wt B(k\theta) \setminus \wt B((k-1)\theta)$ if and only if $x_0 = 0$ and $(x_j = 0 \text{ or } \bar{x}_j = 0)$ for each $j = 1, \dots, n$.
This condition is equivalent to that $b$ does not belong to the left-hand side.
\end{proof}

\begin{lem}\label{lem:dfree}
Let $k$ be a nonnegative integer.
The multiplicity of every element of $\wt B(k\theta) \setminus \wt B((k-1)\theta)$ in $B(k\theta)$ is one. 
\end{lem}

\begin{proof}
Let $b = (x_1, \dots, \bar{x}_1) \in B(k\theta)$ and $\mu = \wt b \in \wt B(k\theta) \setminus \wt B((k-1)\theta)$.
If we set
\[J_+ = \{ j \mid m_j(\mu) > 0\},\]
\[J_- = \{ j \mid m_j(\mu) < 0\},\]
then 
\begin{align*}
x_j &=
\begin{cases}
m_j(\mu) & \text{ for } j \in J_+,\\
0        & \text{ for } j \notin J_+,
\end{cases}\\
\bar{x}_j &=
\begin{cases}
-m_j(\mu) & \text{ for } j \in J_-,\\
0         & \text{ for } j \notin J_-,
\end{cases}\\
x_0 &= 0.
\end{align*}
Hence $b$ is uniquely determined by its weight.
\end{proof}

\begin{thm}\label{thm:dmain}
Let $b \in B(k\theta) \subset B_l$ and assume $b \notin \bigcup_{j=1}^n \Ima \Psi_j$.
\begin{enumerate}
\item We have $\f{0} b=0$ if and only if $k=l$ and $\wt b + \theta \notin \wt B(l\theta)$. 

\item If $\wt b + \theta \in \wt B((k+1)\theta) \setminus \wt B(k\theta)$, then $\f{0} b \in B((k+1)\theta) \subset B_l$.

\item If $\wt b + \theta \in \wt B((k-1)\theta) \setminus \wt B((k-2)\theta)$, then $\f{0} b \in B((k-1)\theta) \subset B_l$.
\end{enumerate}
Moreover, the element $\f{0}b$ is uniquely determined by its weight in each case.
\end{thm}

\begin{proof}
The assertions are immediate from the formula on the action of $\f{0}$, Lemma~\ref{lem:dwt} and Lemma~\ref{lem:dfree}.
\end{proof}

\end{document}